\documentclass{gtart}      
\usepackage{amsmath,amssymb,amscd,rlepsf}  
\input gtmonout
\volumenumber{1}
\volumeyear{1998}
\volumename{The Epstein birthday schrift}
\pagenumbers{383}{411}
\received{21 August 1997}
\published{27 October 1998}
\papernumber{19}

\newcommand{\col}{\hbox{$\,:\,$}}
\renewcommand{\H}{{\mathbb H}}
\renewcommand{\Im}{\operatorname{Im}}
\newcommand{\Ker}{\operatorname{Ker}}
\newcommand{\dehn}{{\delta}}
\newcommand{\dehnker}{\operatorname{\mathcal D}}

\newcommand{\SC}{\operatorname{\mathcal P}}
\newcommand{\SCS}{\operatorname{\mathcal P}(\S^3)}
\newcommand{\SCH}{\operatorname{\mathcal P}(\H^3)}

\newcommand{\CS}{\operatorname{CS}}
\newcommand{\prebloch}{\P}
\newcommand{\bloch}{\B}
\newcommand{\eprebloch}{\widehat\P}
\newcommand{\ebloch}{\widehat\B}
\newcommand{\eeprebloch}{\mathcal E\P}
\newcommand{\eebloch}{\mathcal E\B}

\newcommand{\PSL}{\operatorname{PSL}}
\newcommand{\SL}{\operatorname{SL}}
\newcommand{\SU}{\operatorname{SU}}

\newcommand{\vol}{\operatorname{vol}}

\renewcommand{\P}{{\mathcal P}}
\newcommand{\B}{{\mathcal B}}
\newcommand{\C}{{\mathbb C}}
\newcommand{\Cover}{\widehat\C}
\newcommand{\E}{{\mathbb E}}
\newcommand{\Z}{{\mathbb Z}}
\newcommand{\Q}{{\mathbb Q}}
\newcommand{\R}{{\mathbb R}}
\newcommand\bfw{{\bf w}}
\renewcommand{\S}{{\mathbb S}}
\newcommand{\CP}{{\C\mathbb P}}
\newcommand{\X}{{\mathbb X}}

\newcommand{\calR}{{\mathcal R}}

\newtheorem{theorem}{Theorem}[section]
\newtheorem{theoremq}[theorem]{Theorem?}
\newtheorem{lemma}[theorem]{Lemma}
\newtheorem{proposition}[theorem]{Proposition}
\newtheorem{proponition}[theorem]{Proposition/Definition}
\newtheorem{corollary}[theorem]{Corollary}

\theoremstyle{definition}
\newtheorem{definition}[theorem]{Definition}

\newtheorem{conjecture}[theorem]{Conjecture}

\newtheorem{question}[theorem]{Question}
\newtheorem*{remark*}{Remark}

\begin{document}

\title{Hilbert's 3rd Problem and Invariants of 3--manifolds}
\shorttitle{Hilbert's 3rd problem and invariants of 3--manifolds}
\asciititle{Hilbert's 3rd Problem and Invariants of 3-manifolds}

\author{Walter D Neumann}

\address{Department of Mathematics, The University of
Melbourne\\Parkville, Vic 3052, Australia}

\email{neumann@maths.mu.oz.au}
\begin{abstract}
This paper is an expansion of my lecture for David Epstein's birthday,
which traced a logical progression from ideas of Euclid on subdividing
polygons to some recent research on invariants of hyperbolic
3--manifolds.  This ``logical progression'' makes a good story but
distorts history a bit: the ultimate aims of the characters in the
story were often far from 3--manifold theory.

We start in section 1 with an exposition of the current state of
Hilbert's 3rd problem on scissors congruence for dimension 3.  In
section 2 we explain the relevance to 3--manifold theory and use this
to motivate the Bloch group via a refined ``orientation sensitive''
version of scissors congruence.  This is not the historical motivation
for it, which was to study algebraic $K$--theory of $\C$. Some
analogies involved in this ``orientation sensitive'' scissors
congruence are not perfect and motivate a further refinement in
section \ref{Extended Bloch}. Section \ref{More} ties
together various threads and discusses some questions and conjectures.
\end{abstract}
\asciiabstract{%
This paper is an expansion of my lecture for David Epstein's birthday,
which traced a logical progression from ideas of Euclid on subdividing
polygons to some recent research on invariants of hyperbolic
3-manifolds.  This `logical progression' makes a good story but
distorts history a bit: the ultimate aims of the characters in the
story were often far from 3-manifold theory.

We start in section 1 with an exposition of the current state of
Hilbert's 3rd problem on scissors congruence for dimension 3.  In
section 2 we explain the relevance to 3-manifold theory and use this
to motivate the Bloch group via a refined `orientation sensitive'
version of scissors congruence.  This is not the historical motivation
for it, which was to study algebraic K-theory of C. Some
analogies involved in this `orientation sensitive' scissors
congruence are not perfect and motivate a further refinement in
section 4. Section 5 ties
together various threads and discusses some questions and conjectures.}




\maketitle

\section{Hilbert's 3rd Problem}

It was known to Euclid that two plane polygons of the same area are
related by scissors congruence: one can always cut one of them up into
polygonal pieces that can be re-assembled to give the other.  In the
19th century the analogous result was proved with euclidean geometry
replaced by 2--dimensional hyperbolic geometry or 2--dimensional
spherical geometry.

The 3rd problem in Hilbert's famous 1900 Congress address
\cite{hilbert} posed the analogous question for 3--dimensional
euclidean geometry: are two euclidean polytopes of the same volume
``scissors congruent,'' that is, can one be cut into subpolytopes that
can be re-assembled to give the other.  Hilbert made
clear that he expected a negative answer.  

One reason for the nineteenth century interest in this question was
the interest in a sound foundation for the concepts of area and
volume.  By ``equal area'' Euclid \emph{meant} scissors congruent, and
the attempt in Euclid's Book XII to provide the same approach for
3--dimensional euclidean volume involved what was called an
``exhaustion argument'' --- essentially a continuity assumption ---
that mathematicians of the nineteenth century were uncomfortable with
(by Hilbert's time mostly for aesthetic reasons).

The negative answer that Hilbert expected to his problem was provided
the same year\footnote{In fact, the same answer had been given in 1896
by Bricard, although it was only fully clarified around 1980 that
Bricard was answering an equivalent question --- see Sah's review
85f:52014 (AMS Mathematical Reviews) of \cite{dupont2} for a
concise exposition of this history.} by Max Dehn \cite{dehn}.  Dehn's
answer is delighfully simple in modern terms, so we describe it here
in full.

\begin{definition} Consider the  free $\Z$--module generated by 
the set of congruence classes of 3--dimensional polytopes.  The
{\em scissors congruence group} $\SC(\E^3)$ is the quotient of this
module by the relations of scissors congruence.  That is, if polytopes
$P_1,\dots,P_n$ can be glued along faces to form a polytope $P$ then
we set
$$[P]=[P_1]+\dots+[P_n]\quad\text{in }\SC(\E^3).$$ (A {\em polytope}
is a compact domain in $\E^3$ that is bounded by finitely many planar
polygonal ``faces.'')
\end{definition}

Volume defines a map
$$\vol\co\SC(\E^3)\to \R$$ and Hilbert's problem asks\footnote{%
Strictly speaking this is not quite the same question since two
polytopes $P_1$ and $P_2$ represent the same element of $\SC(\E^3)$ if
and only if they are {\em stably scissors congruent} rather than
scissors congruent, that is, there exists a polytope $Q$ such that
$P_1+Q$ (disjoint union) is scissors congruent to $P_2+Q$.  But, in
fact, stable scissors congruence implies scissors congruence
(\cite{zylev, zylev2}, see \cite{sah-book} for an exposition).} about
injectivity of this map.

Dehn defined a new invariant of scissors congrence, now
called the {\em Dehn invariant}, which can be formulated as a map
$\dehn\co\SC(\E^3)\to\R\otimes\R/\pi\Q$, where the tensor product is
a tensor product of $\Z$--modules (in this case the same as tensor
product as $\Q$--vector spaces).
\begin{definition}
If $E$ is an edge of a polytope $P$
we will denote by $\ell(E)$ and $\theta(E)$ the length of $E$ and
dihedral angle (in radians) at $E$.  For a polytope $P$ we define the
{\em Dehn invariant} $\dehn(P)$ as
$$\dehn(P):=\sum_E\ell(E)\otimes\theta(E)\quad\in\quad
\R\otimes(\R/\pi\Q),\quad\text{sum over all edges $E$ of $P$.}$$
We then extend this linearly to a homomorphism on $\SC(\E^3)$.
\end{definition}

It is an easy but instructive exercise to verify that
\begin{itemize}
\item $\dehn$ is
well-defined on $\SC(\E^3)$, that is, it is compatible with scissors
congruence;
\item $\dehn$ and $\vol$ are independent on $\SC(\E^3)$ in the sense
that their kernels generate $\SC(\E^3)$ (whence
$\Im(\dehn|\Ker(\vol))=\Im(\dehn)$ and
$\Im(\vol|\Ker(\dehn))=\R$);
\item the image of $\dehn$ is uncountable.
\end{itemize}

In particular, $\ker(\vol)$ is not just non-trivial, but even
uncountable, giving a strong answer to Hilbert's question. To give an
explicit example, the regular simplex and cube of equal volume are not
scissors congruent: a regular simplex has non-zero Dehn invariant, and
the Dehn invariant of a cube is zero.

Of course, this answer to Hilbert's problem is really just a start.
It immediately raises other questions:
\begin{itemize}
\item Are volume and Dehn invariant sufficient to classify polytopes
up to scissors congruence?
\item What about other dimensions?
\item What about other geometries?
\end{itemize}

The answer to the first question is ``yes.'' Sydler proved in 1965
that
$$(\vol,\dehn)\co\SC(\E^3)\to \R\oplus(\R\otimes\R/\pi\Q)
$$ 
is injective. Later Jessen \cite{jessen1, jessen2} simplified
his difficult argument somewhat and proved an analogous result for
$\SC(\E^4)$ and the argument has been further simplified in
\cite{dupont-sah-acta}.  Except for these results and the classical
results for dimensions $\le 2$ 
no complete answers are known.  In particular,
fundamental questions remain open about $\SCH$ and $\SCS $.

Note that the definition of Dehn invariant applies with no change to
$\SCH$ and $\SCS $. The Dehn invariant should be
thought of as an ``elementary'' invariant, since it is defined in
terms of 1--dimensional measure.  For this reason (and other reasons
that will become clear later) we are particularly interested in the
kernel of Dehn invariant, so we will abbreviate it: for
$\X=\E^3,\H^3,\S^3$
$$
\dehnker(\X):=\Ker(\dehn\co\SC(\X)\to \R\otimes\R/\pi\Q)
$$

In terms of this notation Sydler's theorem that volume and Dehn
invariant classify scissors congruence for $\E^3$ can be reformulated:
$$\vol\co\dehnker(\E^3)\to\R \quad\text{is injective.}$$
It is believed that volume and Dehn invariant classify
scissors congruence also for hyperbolic and
spherical geometry:
\begin{conjecture}[Dehn Invariant Sufficiency] \label{sufficiency}
$\vol\co\dehnker(\H^3)\to\R$ is injective and 
$\vol\co\dehnker(\S^3)\to\R$ is injective.
\end{conjecture}

On the other hand $\vol\co\dehnker(\E^3)\to\R$ is also surjective,
but this results from the existence of similarity transformations in
euclidean space, which do not exist in hyperbolic or spherical
geometry. In fact, Dupont \cite{dupont1} proved:
\begin{theorem}
$\vol\co\dehnker(\H^3)\to\R$ and 
$\vol\co\dehnker(\S^3)\to\R$ have countable image.
\end{theorem}

Thus the Dehn invariant sufficiency conjecture would imply:
\begin{conjecture}[Scissors Congruence Rigidity]\label{rigidity}
$\dehnker(\H^3)$ and $\dehnker(\S^3)$ are countable.
\end{conjecture}

The following collects results of B\"okstedt, Brun, Dupont, Parry, Sah
and Suslin (\cite{bokstedt-brun-dupont}, \cite{dupont-sah-II},
\cite{sah}, \cite{suslin}).
\begin{theorem}\label{divisibility}
$\SCH$ and $\SCS $ and their subspaces 
$\dehnker(\H^3)$ and $\dehnker(\S^3)$ are uniquely divisible groups,
so they have the structure of\/ $\Q$--vector spaces.  As $\Q$--vector
spaces they have infinite rank.  The rigidity conjecture thus says
$\dehnker(\H^3)$ and $\dehnker(\S^3)$ 
are $\Q$--vector spaces of countably infinite rank.
\end{theorem}

\begin{corollary}\label{image of volume}
The subgroups $\vol(\dehnker(\H^3))$ and\/ $\vol(\dehnker(\S^3))$ of
$\R$ are $\Q$--vector subspaces of countable dimension.
\end{corollary}

\subsection{Further comments}
Many generalizations of Hilbert's problem have been considered, see
eg \cite{sah-book} for an overview. There are generalizations of
Dehn invariant to all dimensions and the analog of the Dehn invariant
sufficiency conjectures have often been made in greater generality,
see eg \cite{sah-book}, \cite{dupont-sah-II}, \cite{goncharov}.
The particular Dehn invariant that we are discussing here is a
codimension 2 Dehn invariant.

Conjecture \ref{sufficiency} appears in various other guises in the
literature.  For example, as we shall see, the $\H^3$ case is
equivalent to a conjecture about rational relations among special
values of the dilogarithm function which includes as a very special
case a conjecture of Milnor \cite{milnor1} about rational linear
relations among values of the dilogarithm at roots of unity.
Conventional wisdom is that even this very special case is a very
difficult conjecture which is unlikely to be resolved in the
forseeable future.  In fact, Dehn invariant sufficiency would imply
the ranks of the vector spaces of volumes in Corollary \ref{image of
volume} are infinite, but at present these ranks are not even proved
to be greater than 1. Even worse: although it is believed that the
volumes in question are always irrational, it is not known if a single
one of them is!

As we describe later, work of Bloch, Dupont, Parry, Sah, Wagoner, and
Suslin connects the Dehn invariant kernels with algebraic $K$--theory
of $\C$, and the above conjectures are then equivalent to standard
conjectures in algebraic $K$--theory.  
In particular, the
scissors congruence rigidity conjectures for $H^3$ and $S^3$ are
together equivalent to the rigidity conjecture for $K_3(\C)$, which
can be formulated that $K_3^{ind}(\C)$ (indecomposable part of
Quillen's $K_3$) is countable.
This conjecture is probably much easier than the Dehn invariant
sufficiency conjecture.

The conjecture about rational relations among special values of the
dilogarithm has been broadly generalized to polylogarithms of all
degrees by Zagier (section 10 of \cite{zagier-conf}).  The connections
between scissors congruence and algebraic $K$--theory have been
generalised to higher dimensions, in part conjecturally, by Goncharov
\cite{goncharov}.

We will return to some of these issues later.  We also refer the
reader to the very attractive exposition in \cite{dupont-sah-new} of
these connnections in dimension 3.

I would like to acknowledge the support of the Australian Research
Council for this research, as well as the the Max--Planck--Institut
f\"ur Mathematik in Bonn, where much of this paper was written.

\section{Hyperbolic 3--manifolds}\label{Hyperbolic 3-manifolds}

Thurston's geometrization conjecture, much of which is proven to be
true, asserts that, up to a certain kind of canonical decomposition,
3--manifolds have geometric structures.  These geometric structures
belong to eight different geometries, but seven of these lead to
manifolds that are describable in terms of surface topology and are
very easily classified. The eighth geometry is hyperbolic geometry
$\H^3$. Thus if one accepts the geometrization conjecture then the
central issue in understanding 3--manifolds is to understand
hyperbolic 3--manifolds.

Suppose therefore that $M=\H^3/\Gamma$ is a hyperbolic 3--manifold.
We will always assume $M$ is oriented and for the moment we will also
assume $M$ is compact, though we will be able to relax this assumption
later.  We can subdivide $M$ into small geodesic tetrahedra, and then
the sum of these tetrahedra represents a class $\beta_0(M)\in \SCH$
which is an invariant of $M$. We call this the {\em scissors
congruence class of $M$}.

Note that when we apply the Dehn invariant to $\beta_0(M)$ the
contributions coming from each edge $E$ of the triangulation sum to
$\ell(E)\otimes2\pi$ which is zero in $\R\otimes\R/\pi\Q$. Thus
\begin{proposition}
The scissors congruence class $\beta_0(M)$ lies in $\dehnker(\H^3)$.\qed
\end{proposition}

How useful is this invariant of $M$?  We can immediately see that it
is non-trivial, since at least it detects volume of $M$:
$$\vol(M)=\vol(\beta_0(M)).$$ Now it was suggested by Thurston in
\cite{thurston3} that the volume of hyperbolic 3--manifolds
should have some close relationship with another geometric invariant,
the Chern--Simons invariant $\CS(M)$. A precise analytic relationship
was then conjectured in \cite{neumann-zagier} and proved in
\cite{yoshida} (a new proof follows from the work discussed here, 
see \cite{neumann1}).  We will not discuss the definition of this
invariant here (it is an invariant of compact riemmanian manifolds,
see
\cite{chern-simons, cheeger-simons}, which was extended also to
non-compact finite volume hyperbolic 3--manifolds by Meyerhoff
\cite{meyerhoff}).  It suffices for the present discussion to know
that for a finite volume hyperbolic 3--manifold $M$ the Chern--Simons
invariant lies in $\R/\pi^2\Z$. Moreover, the combination $\vol(M)+i
\CS(M)\in \C/\pi^2\Z$ turns out to have good analytic properties and
is therefore a natural ``complexification'' of volume for hyperbolic
manifolds.  Given this intimate relationship between volume and
Chern--Simons invariant, it becomes natural to ask if $CS(M)$ is also
detected by $\beta_0(M)$.

The answer, unfortunately, is an easy ``no.''  The point is that
$CS(M)$ is an orientation sensitive invariant: $CS(-M)=-CS(M)$, where
$-M$ means $M$ with reversed orientation.  But, as Gerling pointed out
in a letter to Gauss on 15 April 1844: scissors congruence cannot see
orientation because any polytope is scissors congruent to its mirror
image\footnote{Gauss, Werke, Vol.\ 10, p.\ 242; the argument for a
tetrahedron is to barycentrically subdivide by dropping perpendiculars
from the circumcenter to each of the faces; the resulting $24$
tetrahedra occur in $12$ mirror image pairs.}. Thus
$\beta_0(-M)=\beta_0(M)$ and there is no hope of $CS(M)$ being
computable from $\beta_0(M)$.  This raises the question:

\begin{question}\label{question1} Is there some way to repair 
the orientation insensitivity of scissors congruence and thus capture
Chern--Simons invariant?
\end{question}

The answer to this question is ``yes'' and lies in
the so called ``Bloch group,'' which was invented for entirely
different purposes by Bloch (it was put in final form by Wigner and
Suslin).  To explain this we start with a result of Dupont and Sah
\cite{dupont-sah-II} about ideal polytopes --- hyperbolic polytopes 
whose vertices are at infinity (such polytopes exist in hyperbolic
geometry, and still have finite volume).

\begin{proposition}\label{ideal suffices}
Ideal hyperbolic tetrahedra represent elements in $\SCH$ and,
moreover, $\SCH$ is generated by ideal tetrahedra.
\end{proposition}

To help understand this proposition observe that if $ABCD$ is a
non-ideal tetrahedron and $E$ is the ideal point at which the
extension of edge $AD$ meets infinity then $ABCD$ can be represented
as the difference of the two tetrahedra $ABCE$ and $DBCE$, each of
which have one ideal vertex.  We have thus, in effect, ``pushed'' one
vertex off to infinity.  In the same way one can push a second and
third vertex off to infinity, \dots and the fourth, but this is rather
harder.  Anyway, we will accept this proposition and discuss its
consequence for scissors congruence.

The first consequence is a great gain in convenience: a non-ideal
tetrahedron needs six real parameters satisfying complicated
inequalities to characterise it up to congruence while an ideal
tetrahedron can be neatly characterised by a single complex parameter
in the upper half plane.  

We shall denote the standard compactification of $\H^3$ by $\overline
\H^3 = \H^3\cup\CP^1$. An ideal simplex $\Delta$ with vertices
$z_1,z_2,z_3,z_4\in\CP^1=\C\cup\{\infty\}$ 
is determined up to congruence by the cross-ratio
$$z=[z_1\col z_2\col z_3\col z_4]=\frac{(z_3-z_2)(z_4-z_1)}{(z_3-z_1)(z_4-z_2)}.$$
Permuting the vertices by an even (ie
orientation preserving) permutation replaces $z$ by one of
$$
z,\quad z'=\frac 1{1-z}, \quad\text{or}\quad z''=1-\frac 1z.
$$
The parameter $z$ lies in the upper half plane of $\C$ if the
orientation induced by the given ordering of the vertices agrees with
the orientation of $\H^3$. 

There is another way of describing the cross-ratio parameter $z=
[z_1\col z_2\col z_3\col z_4]$ of a simplex. The group of orientation preserving
isometries of $\H^3$ fixing the points $z_1$ and $z_2$ is isomorphic
to the multiplicative group $\C^*$ of nonzero complex numbers. The
element of this $\C^*$ that takes $z_4$ to $z_3$ is $z$. Thus the
cross-ratio parameter $z$ is associated with the edge $z_1z_2$ of the
simplex.  The parameter associated in this way with the other two
edges $z_1z_4$ and $z_1z_3$ out of $z_1$ are $z'$ and $z''$
respectively, while the edges $z_3z_4$, $z_2z_3$, and $z_2z_4$ have
the same parameters $z$, $z'$, and $z''$ as their opposite edges. See
figure~1.  This description makes clear that the dihedral angles at the
edges of the simplex are $\arg(z)$, $\arg(z')$, $\arg(z'')$
respectively, with opposite edges having the same angle.
\begin{figure}[htbp]
\centerline{
\relabelbox\small
\epsfxsize.35\hsize\epsffile{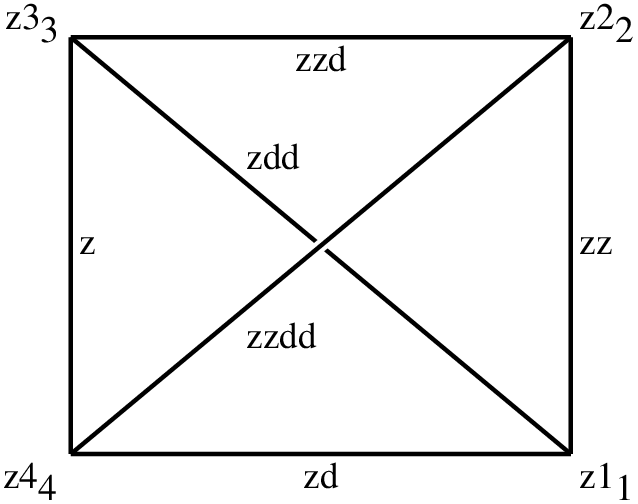}
\relabel{z}{$z$}
\relabel{zz}{$z$}
\relabel{z1}{$z_1$}
\relabel{z2}{$z_2$}
\relabel{z3}{$z_3$}
\adjustrelabel <1pt,0pt> {z4}{$z_4$}
\adjustrelabel <0pt,-2pt> {zd}{$z'$}
\adjustrelabel <0pt,-1.5pt> {zzd}{$z'$}
\relabel{zdd}{$z''$}
\relabel{zzdd}{$z''$}
\endrelabelbox}
\vglue 2mm
\centerline{\small Figure 1}
\end{figure}

Now suppose we have five points
$z_0,z_1,z_2,z_3,z_4\in\CP^1=\C\cup\{\infty\}$.  Any four-tuple of
these five points spans an ideal simplex, and the convex hull of these
five points decomposes in two ways into such simplices, once into two
of them and once into three of them.  We thus get a scissors
congruence relation equating the two simplices with the three
simplices.  It is often called the ``five-term relation.''  To express
it in terms of the cross-ratio parameters it is convenient first to
make an orientation convention.

We allow simplices whose vertex ordering does not agree with the
orientation of $\H^3$ (so the cross-ratio parameter is in the lower
complex half-plane) but make the convention that this represents the
negative element in scissors congruence. An odd permutation of the
vertices of a simplex replaces the cross-ratio parameter $z$ by
$$
\frac 1z, \quad\frac z{z-1},\quad\text{or}\quad 1-z,
$$
so if we denote by $[z]$ the element in $\SCH$ represented by
an ideal simplex with parameter $z$, then our orientation rules say:
\begin{equation}
[z]=[1-\frac 1z]=[\frac 1{1-z}]=-[\frac1z]=-[\frac{z-1}z]=-[1-z].
\label{invsim}
\end{equation}
These orientation rules make the five-term scissors congruence
relation described above particularly easy to state:
$$\sum_{i=0}^4(-1)^i[z_0\col \dots\col \hat{z_i}\col \dots\col z_4]=0.$$ 
The cross-ratio
parameters occuring in this
formula can be expressed in terms of the first two as
\begin{gather*}
{}[z_1\col z_2\col z_3\col z_4] =:x\qquad[z_0\col z_2\col z_3\col z_4]=:y\\
[z_0\col z_1\col z_3\col z_4]=\frac yx\quad
[z_0\col z_1\col z_2\col z_4]=\frac {1-x^{-1}}{1-y^{-1}}\quad
[z_0\col z_1\col z_2\col z_3]=\frac {1-x}{1-y}
\end{gather*}
so the five-term relation can also be written:
\begin{equation}
[x]-[y]+[\frac yx]-[\frac{1-x^{-1}}{1-y^{-1}}]+[\frac{1-x}{1-y}]=0.
\label{5term}
\end{equation}

We lose nothing if we also allow degenerate ideal simplices whose
vertices lie in one plane so the parameter $z$ is real (we always
require that the vertices are distinct, so the parameter is in
$\R-\{0,1\}$), since the five-term relation can be used to express
such a ``flat'' simplex in terms of non-flat ones, and one readily
checks no additional relations result.  Thus we may take the
parameter $z$ of an ideal simplex to lie in $\C-\{0,1\}$ and every
such $z$ corresponds to an ideal simplex.

One can show that relations (\ref{invsim}) follow from the five-term
relation (\ref{5term}), so we consider the quotient
$$\prebloch(\C):=\Z\langle\C-\{0,1\}\rangle/(\text{five-term relations
(\ref{5term}))}$$ of the free $\Z$--module on $\C-\{0,1\}$. Proposition
\ref{ideal suffices} can be restated that there is a
natural surjection $\prebloch(\C)\to\SCH$.  In fact Dupont and Sah
(loc.\ cit.) prove:
\begin{theorem}\label{ideal really suffices}
The scissors congruence group $\SCH$ is the quotient of $\prebloch(\C)$
by the relations $[z]=-[\overline z]$ which identify each ideal
simplex with its mirror image\footnote{The minus sign in this relation
comes from the orientation convention described earlier.}.
\end{theorem}

Thus $\prebloch(\C)$ is a candidate for the orientation sensitive
scissors congruence group that we were seeking. Indeed, it turns out
to do (almost) exactly what we want.

The analog of the Dehn invariant has a particularly elegant expression
in these terms. First note that the above theorem expresses
$\SCH$ as the ``imaginary part'' $\prebloch(\C)^-$ (negative
co-eigenspace under conjugation\footnote{$\prebloch(\C)$ turns out to be a
$\Q$--vector space and is therefore the sum of its $\pm1$ eigenspaces,
so ``co-eigenspace'' is the same as ``eigenspace.''}) of $\prebloch(\C)$.
\begin{proponition}\label{complex dehn} 
The Dehn invariant $\dehn\co\SCH\to\R\otimes\R/\pi\Q$ is twice the
``imaginary part'' of the map
$$\dehn_\C\co\prebloch(\C)\to\C^*\wedge\C^*,\quad [z]\mapsto
(1-z)\wedge z$$ 
so we shall call this map the ``complex Dehn
invariant.''  We denote the kernel of complex Dehn invariant
$$\bloch(\C):=\Ker(\dehn_\C),$$ and call it the ``Bloch group of
$\C$.''
\end{proponition}
(We shall explain this proposition further in an appendix to this
section.)

A hyperbolic 3--manifold $M$ now has an ``orientation sensitive
scissors congruence class'' which lies in this Bloch group and 
captures both volume and Chern--Simons invariant of $M$.  Namely, there
is a map
$$\rho\co \bloch(\C)\to\C/\pi^2\Q$$ introduced by Bloch and Wigner
called the {\em Bloch regulator map}, whose imaginary part is the
volume map on $\bloch(\C)$, and one has:
\begin{theorem}[\cite{neumann-yang3}, \cite{dupont1}]
Let $M$ be a complete oriented hyperbolic 3--manifold of finite
volume.  Then there is a natural class $\beta(M)\in\bloch(\C)$
associated with $M$ and $\rho(\beta(M))=\frac 1i (\vol(M)+i\CS(M))$.
\end{theorem}

This theorem answers Question \ref{question1}.  But there are still
two aesthetic problems:
\begin{itemize}
\item The Bloch regulator $\rho$ plays the r\^ole for orientation
sensitive scissors congruence that volume plays for ordinary scissors
congruence.  But $\vol$ is defined on the whole scissors congruence
group $\SCH$ while $\rho$ is only defined on the kernel $\bloch(\C)$
of complex Dehn invariant.
\item The Chern--Simons invariant $\CS(M)$ is an invariant in $\R/\pi^2\Z$
%\marginpar{coefficient?} 
but the invariant $\rho(\beta(M))$ only computes it in
$\R/\pi^2\Q$.
\end{itemize}
We resolve both these problems in section \ref{Extended Bloch}.

We describe the Bloch regulator map $\rho$ later.  It would be a
little messy to describe at present, although its imaginary part
(volume) has a very nice description in terms of ideal
simplices. Indeed, the volume of an ideal simplex with parameter $z$
is $D_2(z)$, where $D_2$ is the so called ``Bloch--Wigner dilogarithm
function'' given by: $$D_2(z) = \Im \ln_2(z) + \log |z|\arg(1-z),\quad
z\in \C -\{0,1\}$$ and $\ln_2(z)$ is the classical dilogarithm
function. It follows that $D_2(z)$ satisfies a functional equation
corresponding to the five-term relation (see below).

\subsection{Further comments}
To worry about the second ``aesthetic problem'' above could be
considered rather greedy. After all, $CS(M)$ takes values in
$\R/\pi^2\Z$ which is the direct sum of $\Q/\pi^2\Z$ and uncountably
many copies of $\Q$, and we have only lost part of the former summand.
However, it is not even known if the Chern--Simons invariant takes
\emph{any} non-zero values\footnote{According to J Dupont, Jim
Simons deserted mathematics in part because he could not resolve this
issue!} in $\R/\pi^2\Q$. As we shall see, this would be implied by the
sufficiency of Dehn invariant for $\S^3$ (Conjecture
\ref{sufficiency}).

The analogous conjecture in our current situation is:

\begin{conjecture}[Complex Dehn Invariant Sufficiency] 
\label{csufficiency}
$\rho\co\bloch(\C)\to\C/\pi^2\Q$ is injective.
\end{conjecture}

Again, the following is known by work of Bloch:
\begin{theorem}
$\rho\co\bloch(\C)\to\C/\pi^2\Q$ has countable image.
\end{theorem}

Thus the complex Dehn invariant sufficiency conjecture would imply:
\begin{conjecture}[Bloch Rigidity]\label{bloch rigidity}
$\bloch(\C)$ is countable.
\end{conjecture}

\begin{theorem}[\cite{suslin, suslin1}]\label{cdivisibility}
$\prebloch(\C)$ and its subgroup $\bloch(\C)$ are uniquely divisible
groups, so they have the structure of $\Q$--vector spaces.  As
$\Q$--vector spaces they have infinite rank.
\end{theorem}

Note that the Bloch group $\bloch(\C)$ is defined purely algebraically
in terms of $\C$, so we can define a Bloch group $\bloch(k)$
analogously\footnote{Definitions of $\bloch(k)$ in the literature vary
in ways that can mildly affect its torsion if $k$ is not algebraically
closed.} for any field $k$. This group $\bloch(k)$ is uniquely
divisible whenever $k$ contains an algebraically closed field.

It is not hard to see that the rigidity conjecture \ref{bloch
rigidity} is equivalent to the conjecture that
$\bloch(\overline\Q)\to\bloch(\C)$ is an isomorphism (here
$\overline\Q$ is the field of algebraic numbers; it is known that
$\bloch(\overline\Q)\to\bloch(\C)$ is injective). Suslin has
conjectured more generally that $\bloch(k)\to\bloch(K)$ is an
isomorphism if $k$ is the algebraic closure of the prime field in $K$.
Conjecture \ref{csufficiency} has been made in greater generality by
Ramakrishnan \cite{ramakrishnan} in the context of algebraic
$K$--theory.

Conjectures \ref{csufficiency} and \ref{bloch rigidity} are in fact
equivalent to the Dehn invariant sufficiency and rigidity conjectures
\ref{sufficiency} and \ref{rigidity} respectively for $\H^3$ and
$\S^3$ together.  This is because of the following theorem which
connects the various Dehn kernels.  It also describes the connections
with algebraic $K$--theory and homology of the lie group $\SL(2,\C)$
considered as a discrete group. It collates results of of Bloch,
B\"okstedt, Brun, Dupont, Parry and Sah and Wigner (see
\cite{bokstedt-brun-dupont} and \cite{dupont-parry-sah}).
\begin{theorem}
There is a natural exact sequence
$$0\to\Q/\Z\to H_3(\SL(2,\C)) \to \bloch(\C)\to 0.$$ 
Moreover there are natural isomorphisms:

$H_3(\SL(2,\C))\phantom{-}\cong K_3^{ind}(\C)$,

$H_3(\SL(2,\C))^-\cong \bloch(\C)^-\cong \dehnker(\H^3)$,

$H_3(\SL(2,\C))^+\cong \dehnker(\S^3)/\Z$ and
$\bloch(\C)^+\cong\dehnker(S^3)/\Q$, 

\noindent where $\Z\subset\dehnker(S^3)$ is
generated by the class of the 3--sphere and $\Q\subset\dehnker(S^3)$
is the subgroup generated by suspensions of triangles in $\S^2$ with
rational angles. 

The Cheeger--Simons map $c_2\co H_3(\SL(2,\C))\to \C/4\pi^2\Z$ of
\cite{cheeger-simons} induces on the one hand the Bloch regulator map
$\rho\co\bloch(\C)\to\C/\pi^2\Q$ and on the other hand its real and
imaginary parts correspond to the volume maps on $\dehnker(\S^3)/\Z$
and $\dehnker(\H^3)$ via the above isomorphisms.
%\marginpar{see email from Dupont}
\end{theorem}

The isomorphisms of the theorem are proved via isomorphisms
$H_3(\SL(2,\C))^-\cong H_3(\SL(2,R))$ and $H_3(\SL(2,\C))^+ \cong
H_3(\SU(2))$.  We have described the geometry of the isomorphism
$\bloch(\C)^-\cong\dehnker(\H^3)$ in Theorem \ref{ideal really
suffices}. The geometry of the isomorphism
$\bloch(\C)^+\cong\dehnker(S^3)/\Q$ remains rather mysterious.

The exact sequence and first isomorphism in the above theorem are
valid for any algebraically closed field of characteristic $0$.
Thus
Conjecture \ref{bloch rigidity} is also equivalent to each of the four:
\begin{itemize}
\item
Is $K_3^{ind}(\overline\Q)\to K_3^{ind}(\C)$ an isomorphism? Is
$K_3^{ind}(\C)$ countable?
\item
Is $H_3(\SL(2,\overline\Q))\to H_3(\SL(2,\C))$ an isomorphism? Is $
H_3(\SL(2,\C))$ countable?
\end{itemize}

The fact that volume of an ideal simplex is given by the Bloch--Wigner
dilogarithm function $D_2(z)$ clarifies why the $\H^3$ Dehn invariant
sufficiency conjecture \ref{sufficiency} is equivalent to a statement
about rational relations among special values of the dilogarithm
function.  Don Zagier's conjecture about such rational relations,
mentioned earlier, is that any rational linear relation among values
of $D_2$ at algebraic arguments must be a consequence of the relations
$D_2(z)=D_2(\overline z)$ and the five-term functional relation for
$D_2$:
$$
D_2(x)-D_2(y)
+D_2(\frac yx)-D_2(\frac{1-x^{-1}}{1-y^{-1}})+D_2(\frac{1-x}{1-y})=0.
$$
Differently expressed, he conjectures that the volume map is injective
on $\prebloch(\overline \Q)^-$.  If one assumes the scissors
congruence rigidity conjecture for $\H^3$ (that
$\bloch(\overline\Q)^-\cong\bloch(\C)^-$) then the Dehn invariant
sufficiency conjecture for $\H^3$ is just that $D_2$ is injective on
the subgroup $\bloch(\overline\Q)^-\subset\prebloch(\overline\Q)^-$,
so under this assumption Zagier's conjecture is much stronger.
Milnor's conjecture, mentioned earlier, can be formulated that the
values of $D_2(\xi)$, as $\xi$ runs through the primitive $n$-th roots
of unity in the upper half plane, are rationally independent for any
$n$.  This is equivalent to injectivity modulo torsion of the volume
map $D_2$ on $\bloch(k_n)$ for the cyclotomic field $k_n=\Q(e^{2\pi
i/n})$.  For this field $\bloch(k_n)^-=\bloch(k_n)$ modulo
torsion. This is of finite rank but $\prebloch(k_n)^-$ is of infinite
rank, so even when restricted to $k_n$ Zagier's conjecture is much
stronger than Milnor's. Zagier himself has expressed doubt that
Milnor's conjecture can be resolved in the forseeable future.

Conjecture \ref{csufficiency} can be similarly formulated as a
statement about special values of a different dilogarithm function,
the ``Rogers dilogarithm,'' which we will define later.

\subsection{Appendix to section 
\ref{Hyperbolic 3-manifolds}: Dehn invariant of ideal polytopes}
% and Proposition \ref{complex dehn}} 

To define the Dehn invariant of an ideal polytope we first cut off
each ideal vertex by a horoball based at that vertex.  We then have a
polytope with some horospherical faces but with all edges finite.  We
now compute the Dehn invariant using the geodesic edges of this
truncated polytope (that is, only the edges that come from the
original polytope and not those that bound horospherical faces). This
is well defined in that it does not depend on the sizes of the
horoballs we used to truncate our polytope.  (To see this, note that
dihedral angles of the edges incident on an ideal vertex sum to a
multiple of $\pi$, since they are the angles of the horospherical face
created by truncation, which is an euclidean polygon. Changing the
size of the horoball used to truncate these edges thus changes the
Dehn invariant by a multiple of something of the form $l\otimes \pi$,
which is zero in $\R\otimes \R/\pi\Q$.)

Now consider the ideal tetrahedron $\Delta(z)$ with parameter $z$.  We
may position its vertices at $0,1,\infty,z$.  There is a Klein 4--group
of symmetries of this tetrahedron and it is easily verified that it
takes the following horoballs to each other:

\begin{itemize}
\item At $\infty$ the horoball $\{(w,t)\in\C\times\R^+|t\ge a\}$;
\item
at $0$ the horoball of euclidean diameter $|z|/a$;
\item
at $1$ the horoball of euclidean diameter $|1-z|/a$;
\item
at $z$ the horoball of euclidean diameter $|z(z-1)|/a$.
\end{itemize}
After truncation, the vertical edges thus have lengths $2\log
a-\log|z|$, $2\log a-\log|1-z|$, and $2\log a-\log|z(z-1)|$
respectively, and we have earlier said that their angles are $\arg(z),
\arg(1/(1-z)), \arg((z-1)/z)$ respectively.  Thus, adding
contributions, we find that these three edges contribute
$\log|1-z|\otimes\arg(z) - \log |z|\otimes \arg(1-z)$ to the Dehn
invariant.  By symmetry the other three edges contribute the same, so
the Dehn invariant is:
$$\dehn(\Delta(z))=2\bigl(\log|1-z|\otimes\arg(z) - \log|z|\otimes
\arg(1-z)\bigr)\in \R\otimes\R/\pi\Q.$$

\begin{proof}[Proof of Proposition \ref{complex dehn}]
To understand the ``imaginary part'' of
$(1-z)\wedge z \in \C^*\wedge\C^*$ we
use the isomorphism
$$\C^*\to\R\oplus\R/2\pi\Z,\quad z\mapsto \log|z|\oplus\arg z,$$
to represent 
$$\begin{aligned}
\C^*\wedge\C^*&=(\R\oplus\R/2\pi\Z)\wedge(\R\oplus\R/2\pi\Z)\\
&=(\R\wedge\R)\oplus(\R/2\pi\Z
\wedge\R/2\pi\Z)\quad\oplus\quad(\R\otimes\R/2\pi\Z)\\
&=(\R\wedge\R)\oplus(\R/\pi\Q\wedge\R/\pi\Q)
\quad\oplus\quad(\R\otimes\R/\pi\Q),
\end{aligned}
$$ (the equality on the third line is because tensoring over $\Z$ with
a divisible group is effectively the same as tensoring over $\Q$).
Under this isomorphism we have 
$$\begin{aligned} (1-z)\wedge z =
\bigl(\log|1-z|&\wedge\log|z|\oplus\arg(1-z)\wedge\arg
z\bigr)\\&\oplus\quad\bigl(\log|1-z|\otimes\arg z -\log
|z|\otimes\arg(1-z)\bigr),\end{aligned}$$
confirming the Proposition \ref{complex dehn}.
\end{proof}

\section{Computing $\beta(M)$}

The scissors congruence invariant $\beta(M)$ turns out to be a very
computable invariant.  To explain this we must first describe the
``invariant trace field'' or ``field of definition'' of a hyperbolic
3--manifold.  Suppose therefore that $M=\H^3/\Gamma$ is a hyperbolic
manifold, so $\Gamma$ is a discrete subgroup of the orientation
preserving isometry group $\PSL(2,\C)$ of $\H^3$.
\begin{definition}\cite{reid}\qua
The {\em invariant
trace field} of $M$ is the subfield of $\C$ generated over $\Q$ by
the squares of traces of elements of $\Gamma$. We will denote it
$k(M)$ or $k(\Gamma)$. \end{definition}

This field $k(M)$ is an algebraic number field (finite extension of
$\Q$) and is a commensurability invariant, that is, it is unchanged on
passing to finite covers of $M$ (finite index subgroups of $\Gamma$).
Moreover, if $M$ is an arithmetic hyperbolic 3--manifold (that is,
$\Gamma$ is an arithmetic group), then $k(M)$ is the field of
definition of this arithmetic group in the usual sense. See
\cite{reid, neumann-reid1}.

Now if $k$ is an algebraic number field then $\bloch(k)$ is isomorphic
to $\Z^{r_2}\oplus\text{(torsion)}$, where $r_2$ is the number of
conjugate pairs of complex embeddings $k\to\C$ of $k$. Indeed, if
these complex embeddings are $\sigma_1,\dots,\sigma_{r_2}$ then a
reinterpretation of a theorem of Borel \cite{borel} about $K_3(\C)$
says:
\begin{theorem}\label{borel regulator}
The ``Borel regulator map''
$$\bloch(k)\to \R^{r_2}$$ induced on generators of $\prebloch(k)$ by
$[z]\mapsto(\vol[\sigma_1(z)],\dots,\vol[\sigma_{r_2}(z)])$ maps
\hfil\break
$\bloch(k)/\text{(torsion)}$ 
isomorphically onto a full lattice in
$\R^{r_2}$.
\end{theorem}

A corollary of this theorem is that an embedding $\sigma\co k\to\C$
induces an embedding $\bloch(k)\otimes\Q\to\bloch(C)\otimes\Q$. (This
is because the theorem implies that an element of $\bloch(k)$ is
determined modulo torsion by the set of volumes of its Galois
conjugates, which are invariants defined on $\bloch(\C)$.) Moreover,
since $\bloch(\C)$ is a $\Q$--vector space,
$\bloch(\C)\otimes\Q=\bloch(\C)$.

Now if $M$ is a hyperbolic manifold then its invariant trace field
$k(M)$ comes embedded in $\C$ so we get an explicit embedding
$\bloch(k(M))\otimes\Q\to\bloch(\C)$ whose image, which is isomorphic to
$\Q^{r_2}$, we denote by $\bloch(k(M))_\Q$.

\begin{theorem}[\cite{neumann-yang2, neumann-yang3}]
The element $\beta(M)$ lies in the subspace
$\bloch(k(M))_\Q\subset\bloch(\C)$.
\end{theorem}

In fact Neumann and Yang show that $\beta(M)$ is well defined in
$\bloch(K)$ for some explicit multi-quadratic field extension $K$ of
$k(M)$, which implies that $2^c\beta(M)$ is actually well defined in
$\bloch(k(M))$ for some $c$.  Moreover, one can always take $c=0$ if $M$
is non-compact, but we do not know if one can for compact $M$.

In view of this theorem we see that the following data effectively
determines $\beta(M)$ modulo torsion:
\begin{itemize}
\item The invariant trace field $k(M)$.
\item The image of $\beta(M)$ in $\R^{r_2}$ under the Borel regulator
map of Theorem \ref{borel regulator}.
\end{itemize}

To compute $\beta(M)$ we need a collection of ideal simplices that
triangulates $M$ in some fashion.  If $M$ is compact, this clearly
cannot be a triangulation in the usual sense.  In \cite{neumann-yang3}
it is shown that one can use any ``degree one ideal triangulation'' to
compute $\beta(M)$.  This means a finite complex formed of ideal
hyperbolic simplices plus a map of it to $M$ that takes each ideal
simplex locally isometrically to $M$ and is degree one almost
everywhere.  These always exist (see \cite{neumann-yang3} for a
discussion). Special degree one ideal triangulations have been used
extensively in practice, eg in Jeff Weeks' program Snappea
\cite{weeks} for computing with hyperbolic 3--manifolds. Oliver Goodman
has written a program Snap \cite{goodman} (building on Snappea) which
finds degree one ideal triangulations using exact arithmetic in number
fields and computes the invariant trace field and high precision
values for the Borel regulator on $\beta(M)$.

Such calculations can provide numerical evidence for the complex Dehn
invariant sufficiency conjecture. Here is a typical result of such
calculations.
\subsection{Examples}
To ensure that the Bloch group has rank $>1$ we want a field with at
least two complex embeddings.  One of the simplest is the (unique)
quartic field over $\Q$ of discriminant $257$.  This is the field
$k=\Q(x)/(f(x))$ with $f(x)=x^4+x^2-x+1$. This polynomial is
irreducible with roots $\tau_1^{\pm}=0.54742\ldots\pm 0.58565\ldots i$
and $\tau_2^{\pm}=-0.54742\ldots\pm 1.12087\ldots i$.  The field $k$
thus has two complex embeddings $\sigma_1,\sigma_2$ up to complex
conjugation, one with image $\sigma_1(k)=\Q(\tau_1^{-})$ and one with
image $\sigma_2(k)=\Q(\tau_2^-)$.  The Bloch group $\bloch(k)$ is thus
isomorphic to $\Z^2$ modulo torsion.

This field occurs as the invariant trace field 
for two different hyperbolic knot complements in the standard knot
tables up to $8$ crossings, the $6$--crossing knot $6_1$ and the
$7$--crossing knot $7_7$, but the embeddings in $\C$ are different.
For $6_1$ one gets $\sigma_1(k)$ and for $7_7$ one gets $\sigma_2(k)$.
The scissors congruence classes are
$$\begin{aligned} \beta(6_1)&=:\beta_1=2[\frac12(1-\tau^2-\tau^3)] +
[1-\tau] + [\frac12(1-\tau^2+\tau^3)]\in \bloch(k)\\ \beta(7_7)&=:
\beta_2=4[2-\tau-\tau^3]+4[\tau+\tau^2+\tau^3] \in \bloch(k)
\end{aligned}$$
where $\tau$ is the class of $x$ in $k=\Q(x)/(x^4+x^2-x+1)$.  These
map under the Borel regulator $\bloch(k)\to\R^2$ (with respect to the
embeddings $\sigma_1,\sigma_2$) to
$$\begin{aligned} 6_1:\quad& (3.163963228883143983991014716..,
-1.415104897265563340689508587..)\\ 7_7:\quad&
(-1.397088165568881439461453224.., 7.643375172359955478221844448..)
\end{aligned}
$$
In particular, the volumes of these knot complements are 
    $3.1639632288831439..$
and $7.6433751723599554..$ respectively

Snap has access to a large database of small volume compact manifolds.
Searching this database for manifolds whose volumes are small rational
linear combinations of $\vol(\sigma_1(\beta_1))= 3.1639632..$ and
$\vol(\sigma_1(\beta_2))=-1.3970881..$ yielded just eight examples,
three with volume $3.16396322888314..$, four with volume
$4.396672801932495..$ and one with volume $5.629382374981847..$~.  The
complex Dehn invariant sufficiency conjecture predicts (under the
assumption that the rational dependencies found are exact) that these
should all have invariant trace field containing $\sigma_1(k)$.

Checking with Snap confirms that their invariant trace fields equal
$\sigma_1(k)$ and their scissors congruence classes in
$\bloch(k)\otimes\Q$ (computed numerically using the Borel regulator)
are $\beta_1$, $(3/2)\beta_1 + (1/2)\beta_2$, and $2\beta_1+\beta_2$
respectively.

\section{Extended Bloch group}\label{Extended Bloch}

In section \ref{Hyperbolic 3-manifolds} we saw that $\prebloch(\C)$ and
$\bloch(C)$ play a role of ``orientation sensitive'' scissors congruence
group and kernel of Dehn invariant respectively, and that the analog
of the volume map is then the Borel regulator $\rho$.  But, as we
described there, this analogy suffers because $\rho$ is defined on the
Dehn kernel $\bloch(\C)$ rather than on the whole of $\prebloch(\C)$ and
moreover, it takes values in $\C/\pi^2\Q$, rather than in $\C/\pi^2\Z$.

The repair turns out to be to use, instead of $\C-\{0,1\}$, a certain
disconnected $\Z\times\Z$ cover of $\C-\{0,1\}$ to define ``extended
versions'' of the groups $\prebloch(\C)$ and $\bloch(\C)$.  This idea
developed out of a suggestion by Jun Yang.

We shall denote the relevant cover of $\C-\{0,1\}$ by $\Cover$. We
start with two descriptions of it. The second will be a geometric
interpretation in terms of ideal simplices.

Let $P$ be $\C-\{0,1\}$ split along the rays
$(-\infty,0)$ and $(1,\infty)$.  Thus each real number $r$ outside the
interval $[0,1]$ occurs twice in $P$, once in the upper half plane of
$\C$ and once in the lower half plane of $\C$.  We denote these two
occurences of $r$ by $r+0i$ and $r-0i$.  We construct $\Cover$ as an
identification space from $P\times\Z\times\Z$ by identifying
$$
\begin{aligned}
(x+0i, p,q)&\sim (x-0i,p+2,q)\quad\hbox{for each }x\in(-\infty,0)\\
(x+0i, p,q)&\sim (x-0i,p,q+2)\quad\hbox{for each }x\in(1,\infty).\\
\end{aligned}
$$
We will denote the equivalence class of $(z,p,q)$ by $(z;p,q)$.
$\Cover$ has four components:
$$
\Cover=X_{00}\cup X_{01}\cup X_{10}\cup X_{11}
$$
where $X_{\epsilon_0\epsilon_1}$ is the set of $(z;p,q)\in \Cover$ with
$p\equiv\epsilon_0$ and $q\equiv\epsilon_1$ (mod $2$).

We may think of $X_{00}$ as the riemann surface for the function
$\C-\{0,1\}\to\C^2$ defined by $z\mapsto (\log z, -\log (1-z))$.  If
for each $p,q\in\Z$ we take the branch $(\log z + 2p\pi i, -\log (1-z)
+ 2q\pi i)$ of this function on the portion $P\times\{(2p,2q)\}$ of
$X_{00}$ we get an analytic function from $X_{00}$ to $\C^2$.  In the
same way, we may think of $\Cover$ as the riemann surface for the
collection of all branches of the functions $(\log z + p\pi i, -\log
(1-z) + q\pi i)$ on $\C-\{0,1\}$.

We can interpret $\Cover$ in terms of ideal simplices.  Suppose we
have an ideal simplex $\Delta$ with parameter $z\in\C-\{0,1\}$.
Recall that this parameter is associated to an edge of $\Delta$ and
that other edges of $\Delta$ have parameters $$z'=\frac1{1-z},\quad
z''=1-\frac 1z,$$ with opposite edges of $\Delta$ having the same
parameter (see figure~1). Note that $zz'z''=-1$, so the sum $$\log z +
\log z' + \log z''$$ is an odd multiple of $\pi i$, depending on the
branches of $\log$ used.  In fact, if we use the standard branch of
log then this sum is $\pi i$ or $-\pi i$ depending on whether $z$ is
in the upper or lower half plane.  This reflects the fact that the
three dihedral angles of an ideal simplex sum to $\pi$.
\begin{definition}\label{flattening}
We shall call any triple of the form
$$
\bfw=(w_0,w_1,w_2)=(\log z +p\pi i,\log z'+q\pi i,\log z''+r\pi i)
$$ with $$p,q,r\in \Z\quad\text{and}\quad w_0+w_1+w_2=0 $$ a
\emph{combinatorial flattening} for our simplex $\Delta$.  
Thus a combinatorial flattening is an adjustment of each of the three
dihedral angles of $\Delta$ by a multiple of $\pi$ so that the
resulting angle sum is zero.

Each edge $E$ of $\Delta$ is assigned one of the
components $w_i$ of $\bfw$, with opposite edges being assigned the
same component. We call
$w_i$ the \emph{log-parameter} for the edge $E$ and denote it
$l_E(\Delta,\bfw)$. 
\end{definition}

For $(z;p,q)\in \Cover$ we define
$$\ell(z;p,q):=(\log z + p\pi i, -\log(1-z) + q\pi i,
\log(1-z)-\log z - (p+q)\pi i),
$$ 
and $\ell$ is then a map
of $\Cover$ to the set of combinatorial flattenings of simplices.

\begin{lemma}\label{Cover is flattenings} 
This map $\ell$ is a bijection, so $\Cover$ may be identified with
the set of all combinatorial flattenings of ideal tetrahedra.
\end{lemma}

\begin{proof}
We must show that $(w_0,w_1,w_2)= \ell(z;p,q)$ determines $(z;p,q)$.
It clearly suffices to recover $z$. But up to sign $z$ equals
$e^{w_0}$ and $1-z$ equals $e^{-w_1}$, and the knowledge of both $z$
and $1-z$ up to sign determines $z$.
\end{proof}

\subsection{The extended groups}
We shall define a group $\eprebloch(\C)$ as
$\Z\langle\Cover\rangle/$(relations), where the relations in question
are a lift of the five-term relations (\ref{5term}) that define
$\prebloch(\C)$, plus an extra relation that just eliminates an
element of order $2$.

We first recall the situation of the five-term relation (\ref{5term}).
If $z_0,\ldots,z_4$
are five distinct points of $\partial\overline\H^3$, then each choice
of four of five points $z_0,\dots,z_4$ gives an ideal simplex. We
denote the simplex which omits vertex $z_i$ by $\Delta_i$. 
We denote the cross-ratio parameters of these simplices by 
$x_i=[z_0\col \ldots\col \hat{z_i}\col \ldots\col z_4]$. Recall that $(x_0,\dots,x_4)$
can be written in terms of $x=x_0$ and $y=x_1$ as
$$(x_0,\dots,x_4)=\biggl(x,y,\frac yx,\frac{1-x^{-1}}{1-y^{-1}},
\frac{1-x}{1-y}\biggr)$$
The five-term relation was $\sum_{i=0}^4(-1)^i[x_i]=0$, so the lifted
five-term relation will have the form \def\term#1{(x_#1;p_#1,q_#1)}
\begin{equation}\label{5-term}
\sum_{i=0}^4(-1)^i \term i=0
\end{equation}
with certain relations on the $p_i$ and $q_i$.  We need to describe
these relations.

Using the map of Lemma
\ref{Cover is flattenings}, each summand in this relation
(\ref{5-term}) represents a choice $\ell\term i$ of combinatorial
flattening for one of the five ideal simplices. For each edge $E$
connecting two of the points $z_i$ we get a corresponding linear
combination
\begin{equation}\label{edge sums}
\sum_{i=0}^4(-1)^il_E(\Delta_i,\ell\term i)
\end{equation}
of log-parameters (Definition \ref{flattening}), where we put
$l_E(\Delta_i,\ell\term i)=0$ if the line $E$ is not an edge of
$\Delta_i$.  This linear combination has just three non-zero terms
corresponding to the three simplices that meet at the edge $E$. One
easily checks that the real part is zero and the imaginary part can be
interpreted (with care about orientations) as the sum of the
``adjusted angles'' of the three flattened simplices meeting at $E$.

\begin{definition}\label{geometric five term}
We say that the $(x_i;p_i,q_i)$ satisfy the \emph{flattening
condition} if each of the above linear combinations (\ref{edge sums})
of log-parameters is equal to zero. That is, the adjusted angle sum of
the three simplices meeting at each edge is zero.  In this case
relation (\ref{5-term}) is an instance of the \emph{lifted five-term
relation}.
\end{definition}

There are ten edges in question, so the flattening conditions are ten
linear relations on the ten integers $p_i,q_i$.  But these equations
turn out to be linearly dependant, and the space of solutions is
$5$--dimensional.  For example, if the five parameters $x_0,\dots,x_4$
are all in the upper half-plane (one can check that this means $y$ is
in the upper half-plane and $x$ is inside the triangle with vertices
$0,1,y$) then the conditions are equivalent to:
\begin{gather}
p_2=p_1-p_0,\quad p_3=p_1-p_0+q_1-q_0,\quad p_4=q_1-q_0\nonumber\\
q_3=q_2-q_1,\quad q_4 = q_2-q_1-p_0\nonumber \end{gather}
 which express $p_2$, $p_3$,
$p_4$, $q_3$, $q_4$ in terms of $p_0$, $p_1$, $q_0$, $q_1$, $q_2$.
Thus, in this case the lifted five-term relation becomes:
\begin{gather}
(x_0;p_0,q_0)-(x_1;p_1,q_1)+(x_2;p_1-p_0,q_2)-{}\nonumber\\ 
{}-(x_3;p_1-p_0+q_1-q_0,q_2-q_1)+(x_4;q_1-q_0,q_2-q_1-p_0)=0\nonumber
\end{gather}
This situation corresponds to the configuration of figure~2
for the ideal vertices $z_0,\dots,z_4$,
\begin{figure}[hbtp]
\centerline{\relabelbox\small
\epsfxsize.5\hsize\epsffile{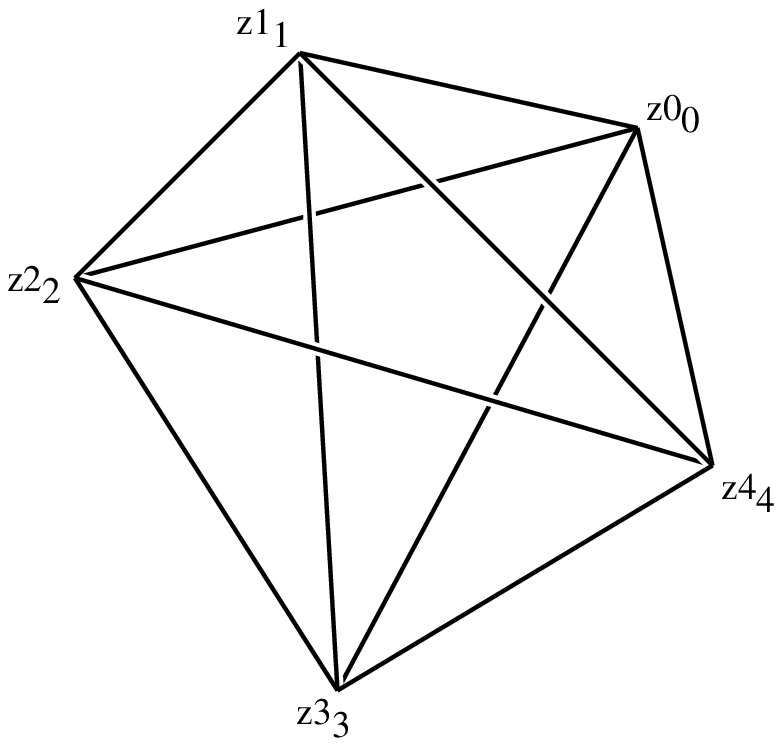}
\relabel{z0}{$z_0$}
\adjustrelabel <9pt, 0pt> {z1}{$z_1$}
\adjustrelabel <4pt, 0pt> {z2}{$z_2$}
\adjustrelabel <5pt, 0pt> {z3}{$z_3$}
\adjustrelabel <0pt, 5pt> {z4}{$z_4$}
\endrelabelbox}
\vglue -7mm
\centerline{\small Figure 2}
\end{figure}
with $z_1$ and $z_3$ on opposite sides of the plane of the triangle
$z_0z_2z_4$ and the line from $z_1$ to $z_3$ passing through the
interior of this triangle.  

\begin{definition}
The extended pre-Bloch group $\eprebloch(\C)$ is the group
$$\eprebloch(\C) := \Z\langle\Cover\rangle/(\text{lifted five-term
relations and the following relation})$$
\begin{equation}\label{transfer}
[x;p,q]+[x;p',q']=[x;p,q']+[x;p',q].
\end{equation}
\end{definition}
(We call relation (\ref{transfer}) the \emph{transfer relation}; one
can show that if one omits it then $\eprebloch(\C)$ is replaced by
$\eprebloch(\C) \oplus\Z/2$, where the $\Z/2$ is generated by the
element $\kappa:=[x,1,1]+[x,0,0]-[x,1,0]-[x,0,1]$, which is
independant of $x$.)

The relations we are using are remarkably natural.  To explain this we
need a beautiful version of the dilogarithm function called the
\emph{Rogers dilogarithm}:
$${\calR}(z)=-\frac12\biggl(\int_0^z\bigl(\frac{\log
t}{1-t}+\frac{\log(1-t)}t\bigr)dt\biggr)-\frac{\pi^2}6.$$ The extra
$-\pi^2/6$ is not always included in the definition but it improves
the functional equation. $\calR(z)$ is singular at $0$ and $1$ and is
not well defined on $\C-\{0,1\}$, but it lifts to an analytic function
\begin{gather}
R:\Cover\to\C/\pi^2\Z\nonumber\\ R(z;p,q)={\calR}(z)+\frac{\pi
i}2(p\log(1-z)+q\log z)  .\nonumber\end{gather} We also consider the map
$$\hat\delta\co\Cover\to\C\wedge\C,\quad \hat\delta(z;p,q)=
\bigl(\log z + p\pi i\bigr)\wedge \bigl(-\log(1-z) +q\pi i\bigr).$$
Relation (\ref{transfer}) is clearly a functional equation for both
$R$ and $\hat\delta$.  It turns out that the same is true for the
lifted five-term relation.  In fact:
\begin{proposition}
If $(x_i;p_i,q_i)$, $i=0,\dots,4$ satisfy the flattening condition, so
$$\sum_{i=0}^4(-1)^i \term i=0$$
is an instance of the lifted five-term relation, then
$$\sum_{i=0}^4(-1)^i R\term i=0$$
and
$$\sum_{i=0}^4(-1)^i \hat\delta\term i=0.$$ Moreover, each of these
equations also characterises the flattening condition.
\end{proposition}

Thus the flattening condition can be defined either geometrically, as
we introduced it, or as the condition that makes the five-term
functional equation for either $R$ or $\hat\delta$ valid.  In any
case, we now have:
\begin{theorem}
$R$ and $\hat\delta$ define maps
$$\begin{aligned}
R\co &\eprebloch(\C)\to \C/\pi^2\Z\\
\hat\delta\co &\eprebloch(\C) \to \C\wedge\C.
  \end{aligned}$$
\end{theorem}

We call $\hat\delta$ the \emph{extended Dehn invariant} and call its
kernel $$\ebloch(\C):=\ker(\hat\delta)$$ the \emph{extended Bloch
group}.  The final step in our path from Hilbert's 3rd problem to
invariants of 3--manifolds is given by the following theorem.

\begin{theorem}\label{final step}
A hyperbolic 3--manifold $M$ has a natural class
$\hat\beta(M)\in\ebloch(\C)$.  Moreover%
%\footnote{At the time of
%writing this formula is only proved modulo $\pi^2/6$, since it
%depends on computations in \cite{dupont-kamber}, where best
%denominators were not sought after.}
, $R(\hat\beta(M))=\frac 1i
(\vol(M)+i\CS(M))\in\C/\pi^2\Z$.
\end{theorem}

To define the class $\hat\beta(M)$ directly from an ideal
triangulation one needs to use a more restrictive type of ideal
triangulation than the degree one ideal triangulations that suffice
for $\beta(M)$. For instance, the triangulations constructed by
Epstein and Penner \cite{epstein-penner} in the non-compact case and
by Thurston \cite{thurston2} in the compact case are of the
appropriate type.  One then chooses flattenings of the ideal simplices
of $K$ so that the whole complex $K$ satisfies certain ``flatness''
conditions. The sum of the flattened ideal simplices then represents
$\hat\beta(M)$ up to a $\Z/6$ correction. The main part of the
flatness conditions on $K$ are the conditions that adjusted angles
around each edge of $K$ sum to zero together with similar conditions
on homology classes at the cusps of $M$. If one just requires these
conditions one obtains $\hat\beta(M)$ up to $12$--torsion.  Additional
mod $2$ flatness conditions on homology classes determine
$\hat\beta(M)$ modulo $6$--torsion.  The final $\Z/6$ correction is
eliminated by appropriately ordering the vertices of the simplices of
$K$.  It takes some work to see that all these conditions can be
satisfied and that the resulting element of $\ebloch(\C)$ is well
defined, see
\cite{neumann, neumann1}.

\section{Comments and questions}\label{More}

\subsection{Relation with the non-extended Bloch group}

What really underlies the above Theorem \ref{final step} is the
\begin{theoremq}\label{theoremq}
There is a natural short exact sequence
$$0\to\Z/2\to H_3(\PSL(2,\C);\Z)\to \ebloch(\C)\to 0.$$\end{theoremq}

The reason for the question mark is that, at the time of writing, the
proof that the kernel is exactly $\Z/2$ has not yet been written down
carefully.

The relationship of our extended groups with the
``classical'' ones is as follows.
\begin{theorem} There is a commutative diagram with exact rows and
columns:
$$
\begin{CD}
&& 0 &&  0 &&  0 \\
&& @VVV @VVV @VVV \\
0 @>>> \mu^* @>>> \C^* @>>> \C^*/\mu^* @>>> 0\\
&& @V\chi|\mu^* VV @V\chi VV @V\xi VV 
%\\
@VVV\\
0 @>>> \ebloch(\C) @>>> \eprebloch(\C) @>\hat\delta>> \C\wedge\C 
%\\
@>>> K_2(\C) @>>> 0\\
&& @VVV @VVV @V\epsilon VV 
%\\
@V=VV\\
0 @>>> \bloch(\C) @>>> \prebloch(\C) @>\delta>> \C^*\wedge\C^* 
%\\
@>>> K_2(\C) @>>> 0\\
&& @VVV @VVV @VVV 
%\\
@VVV\\
&& 0 &&  0 &&  0 
 && 0
\end{CD}
$$
Here $\mu^*$ is the group of roots of unity and the labelled maps
that have not yet been defined are:
$$\begin{aligned}
\chi(z)&=[z,0,1]-[z,0,0]\in
\eprebloch(\C);\\
\xi[z]&=\log z \wedge \pi i;\\
\epsilon(w_1\wedge w_2)&= (e^{w_1}\wedge e^{w_2}).
  \end{aligned}
$$%\marginpar{Check Sign, also for delta}
%(The cokernels of the maps $\delta$ and $\hat\delta$ are $K_2(\C)$.) 
\end{theorem}

%\section{Comments and Questions}\label{More}

\subsection{Extended extended Bloch}
The use of the disconnected cover $\Cover$ of $\C-\{0,1\}$ rather than
the universal abelian cover (the component $X_{00}$ of $\Cover$) in
defining the extended Bloch group may seem unnatural.  If one uses
$X_{00}$ instead of $\Cover$ one obtains extended Bloch groups
$\eeprebloch(\C)$ and $\eebloch(\C)$ which are non-trivial $\Z/2$
extensions of $\eprebloch(\C)$ and $\ebloch(\C)$. Theorem
\ref{theoremq} then implies a natural
\emph{isomorphism} $H_3(\PSL(2,\C);\Z)\to\eebloch(\C)$.  The 
homomorphism of Theorem \ref{theoremq} is given explicitely by
``flattening'' homology classes in the way sketched after Theorem
\ref{final step}, and the isomorphism
$H_3(\PSL(2,\C);\Z)\to\eebloch(\C)$ presumably has a similar explicit
description using ``$X_{00}$--flattenings,'' but we have not yet proved
that these always exist (note that an $X_{00}$--flattening
of a simplex presupposes a choice of a pair of opposite edges of the
simplex; changing this choice turns it into a $X_{01}$-- or
$X_{10}$--flattening).

For the same reason, we do not yet have a simplicial
description of the class $\hat\beta(M)\in\eebloch{\C}$ for a closed
hyperbolic manifold $M$, although this class exists for homological
reasons.  It is essential here that $M$ be closed --- the class
$\hat\beta(M)\in\ebloch(\C)$ almost certainly has no natural lift to
$\eebloch(\C)$ in the non-compact case. 

The Rogers dilogarithm induces a natural map
$R\co\eebloch(\C)\to\C/2\pi^2\Z$, and this is the Cheeger--Simons
class $H_3(PSL(2,\C)\to\C/2\pi^2\Z$ via the above isomorphism.

\subsection{Computing Chern--Simons invariant}
The formula of \cite{neumann} for $CS(M)$ used in the programs Snappea
and Snap uses ideal triangulations that arise in Dehn surgery.  These
triangulations are not of the type mentioned after Theorem \ref{final
step}, but by modifying them one can put them in the desired form and
use Theorem \ref{final step} to compute $\hat\beta(M)$, reconfirming
the formula of \cite{neumann}. The formula computes $CS(M)$ up to a
constant for the infinite class of manifolds that arise by Dehn
surgery on a given manifold.  It was conjectured in \cite{neumann}
that this constant is always a multiple of $\pi^2/6$, and this too is
confirmed. The theorem also gives an independent proof of the relation
of volume and Chern--Simons invariant conjectured in
\cite{neumann-zagier} and proved in
\cite{yoshida}, from which a formula for eta-invariant was also 
deduced in \cite{neumann-meyerhoff} and \cite{ouyang}.

\subsection{Realizing elements in the Bloch group and Gromov norm}

One way to prove the Bloch group rigidity conjecture \ref{bloch
rigidity} would be to show that $\bloch(\C)$ is generated by the
classes $\beta(M)$ of 3--manifolds.  This question is presumably much
harder than the rigidity conjecture, although modifications of it have
been used in attempts on it. More specifically, one can ask 
\begin{question}For which number fields $k$ is
$\bloch(k)_\Q$ generated as a $\Q$ vector space by classes
$\beta(M)$ of 3--manifolds with invariant trace field contained in $k$?
\end{question}

For totally real number fields (ie $r_2=0$) the answer is trivially
``yes'' while for number fields with $r_2=1$ the existence of
arithmetic manifolds again shows the answer is ``yes.'' But beyond
this little is known. In fact it is not even known if for every
non-real number field $k\subset\C$ a 3--manifold exists with
invariant trace field $k$.  (For a few cases, eg multi-quadratic
extensions of $\Q$, the author and A Reid have unpublished
constructions to show the answer is ``yes.'')

Jun Yang has pointed out that ``Gromov norm'' gives an obstruction to
a class $\alpha\in\bloch(\C)$ being realizable as $\beta(M)$
(essentially the same observation also occurs in \cite{reznikov}).  We
define the \emph{Gromov norm} $\nu(\alpha)$ as
$$\nu(\alpha)=\inf\bigl\{\sum\bigl|\frac{n_i}k\bigr|:k\alpha=
\sum n_i[z_i],\quad z_i\in
\C\bigr\},$$ and it is essentially a result of Gromov, with proof given
in \cite{thurston1}, that:
\begin{theorem}
$$|\vol(\alpha)|\le V\nu(\alpha),$$
where $V=1.00149416...$ is the volume of a regular ideal tetrahedron.
If $\alpha=\beta(M)$ for some 3--manifold $M$ then
$$\vol(\alpha)=V\nu(\alpha).$$
\end{theorem}

In particular, since $\nu(\alpha)$ is invariant under the action of
Galois, for $\alpha=\beta(M)$ one sees that the $\vol(M)$ component of
the Borel regulator is the largest in absolute value and equals
$V\nu(\alpha)$.  This suggests the question:
\begin{question}
Is it true for any number field $k$ and for any $\alpha\in\bloch(k)$ that
$V\nu(\alpha)$ equals the largest absolute value of a component of
the Borel regulator of $\alpha$?
\end{question}

This question is rather naive, and at this point we have no evidence
for or against.  Another naive question is the following. For
$\alpha\in\bloch(k)_\Q$, where $k$ is a number field, we can define a
stricter version of Gromov norm by
$$\nu_k(\alpha)=\inf\bigl\{\sum\bigl|\frac{n_i}k\bigr|:
k\alpha=\sum n_i[z_i],\quad z_i\in k\bigr\}.$$
\begin{question}
Is $\nu_k(\alpha)=\nu(\alpha)$ for $\alpha\in\bloch(k)_\Q$?
\end{question}
If not, then $\nu_k$ gives a sharper obstruction to realizing $\alpha$
as $\beta(M)$ since it is easy to show that for $\alpha=\beta(M)$ one
has $\vol(\alpha)=V\nu_K(\alpha)$ for some at most quadratic extension
$K$ of $k$.

\bigskip
{\small \parskip 0pt \leftskip 0pt \rightskip 0pt plus 1fil \def\\{\par}
\sl\theaddress\par
\medskip
\rm Email:\stdspace\tt\theemail\par}
\recd

\end{document}